\documentclass{amsart}
\usepackage{graphics, graphicx}
\usepackage{amsthm, amsmath, amssymb}

\newcommand{\ds}{\displaystyle}

\def\NW{{\sf NW}}\def\NN{{\sf N}}
\def\SE{{\sf SE}}\def\EE{{\sf E}}
\def\SS{{\sf S}}\def\NE{{\sf NE}}
\def\WW{{\sf W}}\def\SW{{\sf SW}}
\def\LY{\mathcal{L}(\Y)}
\def\AA{\mathcal{A}}
\def\BB{\mathcal{B}}

\newcommand\Z{\ensuremath{\mathbb Z}}

\def\Y{\ensuremath{\mathcal Y}} % generic walk set

\def\rev{\mathsf{rev}}
\def\refl{\mathsf{reflect}}
\def\y{Y}

\def\xn{\underline{x}}
\newtheorem{theorem}{Theorem}[section]
\newtheorem{proposition}{Proposition}[section]
\newtheorem{lemma}[theorem]{Lemma}
\newtheorem{conjecture}{Conjecture}
\theoremstyle{definition}

\theoremstyle{remark}

\numberwithin{equation}{section}

\numberwithin{equation}{section}
\newtheoremstyle{marnidef}{\topsep}{\topsep}%
     {}%         Body font
     {}%         Indent amount (empty = no indent, \parindent = para indent)
     {\bfseries}% Thm head font
     {.  }%        Punctuation after thm head
     {.5em}%     Space after thm head (\newline = linebreak)
     {\thmname{#1}\thmnumber{#2 }\thmnote{{\em #3}}}%         Thm head spec
\theoremstyle{marnidef}
\newtheorem{defn}{Definition }[section]

\author{Marni Mishna}
\address{Dept. of Mathematics, Simon Fraser University, 8888 University Dr. Burnaby, Canada}
\email{mmishna@sfu.ca}
\urladdr{http://www.math.sfu.ca/~mmishna}
\title{Classifying lattice walks restricted to the quarter plane}
\subjclass[2000]{Primary 82B41; Secondary 05A15}
\keywords{enumeration, lattice walks, D-finite functions, holonomic functions}
\begin{document}
\begin{abstract}
  This work considers lattice walks restricted to the quarter plane, with
  steps taken from a set of cardinality three. We present a complete classification of
  the generating functions of these walks with respect to the classes
  algebraic, transcendental holonomic and non-holonomic. The principal
  results are a new algebraic class related to Kreweras' walks; two
  new non-holonomic classes; and enumerative data on some other
  classes. These results provide strong evidence for conjectures which
  use combinatorial criteria to classify the generating functions all
  nearest neighbour walks in the quarter plane.
\end{abstract}

\maketitle

\section*{Introduction}
The interest in an enumerative approach to lattice walks under various
different types of restrictions has risen recently, (see \cite{MBM02,
  MBM05, BoPe03, JaPrRe06}) and there is increased need for global
approaches and results. A few such studies have been performed. For
example, in the case of the one dimensional lattice walks, Flajolet and
Banderier~\cite{BaFl02} examine the nature of their generating
functions, and provide general results of asymptotic analysis. For two
dimensional walks in the quarter plane, we find a collection of case
analyses~ \cite{MBM05, BoPe03, Gessel86, Regev81}, and the goal here
is to try to determine more general characterizations of these walks,
based on the nature of their generating functions. Essentially we are
interested to know which walks have holonomic generating functions,
that is, when does the generating function satisfy
systems of independent linear differential equations with polynomial
coefficients. The answer has important repercussions for
sequence generation, asymptotics, amongst other enumerative questions.

Unfortunately we do not completely succeed in giving a 
characterization of all walks, but we do uncover some interesting
patterns, and present variety of applications of the kernel method. We
derive enumerative information to support a new conjecture
on the combinatorial conditions required for a class of walks to
possess holonomic (or, D-finte) generating functions.

\section{Walks and their generating functions} 
\subsection{Next nearest neighbour walks}
The walks of interest here are known as {\em next nearest neighbour
walks}. Precisely, they use movements on the integer lattice where
each step is from some fixed set $\Y\subseteq\{\pm1,
0\}^2\setminus\{(0,0)\}$, which we also specify by the compass
directions~$\{\NN, \NE, ..., \WW, \NW\}$. Such a set $\Y$ is called a
{\em step set}. Here we shall consider nearest neighbour walks
exclusively, unless explicitly mentioned otherwise. A {\em walk in the
quarter plane} is a sequence of steps~$w$ in $\Y^*$, $w=w_1, w_2, \ldots, w_n$,
such that for each $k\leq n$, the vector sum $(x_k,y_k)=\sum_{i=1}^k w_i$
satisfies $x_k\geq 0, y_k\geq 0$, that is, it remains in the first
quadrant. We shall denote the set of all valid walks with steps
from~$\Y$ by~$\LY$.  We can consider this as a formal language
(in the sense of theoretical computer science) over the alphabet $\Y$
with the horizontal and vertical conditions as prefix conditions on
any word in the language.

\subsection{Complete and counting generating functions}
Fix some step set~$\Y$. We associate to~$\LY$ two power series:~$W_\Y(t)$ a
counting (ordinary, univariate) generating function and a complete
(multivariate) generating function. The series~$W_\Y(t)$, is a
formal power series where the coefficient of $t^n$ is the
number of walks of length~$n$. The complete generating function
$Q_\Y(x,y;t)$ encodes more information. The coefficient of $x^iy^jt^n$
in~$Q_\Y(x,y;t)$ is the number of walks of length $n$ ending at
the point~$(i,j)$. Remark that the specialization $x=y=1$ in the complete
generating function is precisely $W_\Y(t)$,
i.e. $Q_\Y(1,1;t)=W_\Y(t)$. In both cases, when $\Y$ is clear it is often dropped as an index.

In part, our interest in the complete generating function stems from the
fact that if it is in a particular functional class, then generally so
is $W_\Y$. Furthermore, we can determine a useful functional equation that it
satisfies. The {\em fundamental equation\/} satisfied by a complete
generating function is determined from the recursive definition that a
walk of length~$n$ is a walk of length $n-1$ plus a step. For a walk
ending on the $x$- or $y$-axis, it is possible that not all of the
directions from $\Y$ will be permissible for the next step, and thus we
subtract out these sub-series. 
\begin{defn}[Fundamental equation] The fundamental equation of the
complete generating function of walks with step set $\Y$ is given
by\footnote{$\chi[P]=1$ if $P$ is true and $0$ otherwise; $\bar
x=\frac1x$; $\bar y=\frac1y$} 
\begin{align}\label{eqn:fund}
 Q(x,y;t)=1&+
             \sum_{(i, j)\in\Y}\! tx^i y^j Q(x,y;t)
            -t\bar y\hspace{-3mm}\sum_{i:(i, -1)\in\Y}\!x^i Q(x,0;t)  
            -t\bar x\hspace{-3mm}\sum_{j:(-1, j)\in\Y}\!y^j  Q(0,y;t)\\
          &\nonumber +\chi [(-1,-1)\in\Y]\,t\bar x\bar yQ(0,0;t).
\end{align}
\end{defn}
\subsection{Classifying formal series} 
These two series fall into the following classes of functions. Let
$\xn=x_1, x_2, \ldots, x_n$.  A multivariate generating
function~$G(\xn)$ is {\em algebraic\/} if there exists a multivariate
polynomial $P(\xn,y)$ such that $P(\xn, G(\xn))=0$.
Flajolet~\cite{Flajolet87} summarizes a number of criteria which imply
the transcendence of a series. One which we shall use here is a
consequence of his Theorem D: If the Taylor coefficient of $z^n$ of a
function $f(z)$ (analytic at the origin) is asymptotically
equivalent to $\gamma \beta^n n^r$, and further if  $r$ is
irrational or a negative integer; or if $\beta$ is transcendental; or if
$\gamma\Gamma(r+1)$ transcendental, then $f(z)$ is transcendental.

\begin{defn}[holonomic function]
A multivariate series $G(\xn)$ is {\em holonomic\/} if the
vector space generated by its partial derivatives (and their
iterates), over rational functions of $\xn$ is finite
dimensional. This is equivalent to the existence of $n$ partial
differential equations of the form \[0=p_{0, i} f(\xn)+ p_{1,i} \frac{\partial
f(\xn)}{\partial x_i}+\ldots + p_{d_i, i} \frac{\partial^{d_i}
f(\xn)}{(\partial x_i)^{d_i}} ,\] for $i$ satisfying $1\leq i\leq n$,
and where the $p_{j, i}$ are all polynomials in $\xn$. Holonomic functions
are also known as {\em D-finite\/} functions.
\end{defn}

Algebraic functions are always holonomic, but as $\exp(x)$ is
holonomic and  transcendental, we see this containment is strict. The
closure properties of these two classes are presented by
Stanley~\cite[Ch. 6]{Stanley99}.  In particular, here we shall use the
fact that if $F(\xn)$ is holonomic with respect to the $x_i$, and if
the algebraic substitution $x_i=y_i(z_1, \ldots, z_k)\equiv y_i$ makes
sense as a power series substitution, then $f(y_1, \ldots, y_n)$ is
holonomic with respect to the $z_i$ variables.

The goal of this work is to determine when $W_\Y$ and $Q_\Y$ are holonomic,
algebraic, or neither. The next two theorems are a model example of the kind of
result we aspire to emulate, from the classification point of view.

\begin{theorem}[Half-plane condition~\cite{BaFl02}]
\label{thm:hpc}
Let~$\Y$ be a subset of $\{\pm 1,0\}^2$. The complete generating
series~$Q(x,y;t)$ for walks that start at $(0,0)$, take their steps
in~$\Y$ and stay in a half plane is algebraic. 
\end{theorem}

A step set $\Y$ is {\em $y$-axis symmetric\/} (resp. $x$-axis
symmetric) if~$(i,j)\in \Y$ implies that
$(i, -j)\in \Y$ (resp. $(-i, j)\in\Y$). 
\begin{theorem}[Bousquet-M\'elou, Petkov\v sek~\cite{ MBM03, MBM05}]
\label{thm:sym}
 Let $\Y$ be a
finite subset of $\{\pm 1, 0\}\times \Z\setminus\{(0,0)\}$ that is
symmetric with respect to the $y$-axis. Then the complete generating
function~$Q(x, y; t)$ for walks that start from $(0,0)$, take their
steps in $\Y$, and stay in the first quadrant is holonomic.
\end{theorem}

For example, $\Y=\{\EE, \NW, \SW\}$ satisfies the conditions of Theorem~\ref{thm:sym}, and
thus $G_\Y(x, y;t)$ is a holonomic function. Naturally, an analogous
result for step sets from $\Z\times\{\pm 1, 0\}\setminus\{(0,0)\}$
which are $x$-axis symmetric is also true.

\subsection{Combinatorial operations and conditions} Here we introduce two
combinatorial operations which act on the step sets and which play an
important role. First, we have $\refl$, which switches the
coordinates, $\refl(x,y)=(y,x)$, effectively flipping each step across
the line~$x=y$. This operator preserves both algebraicity and
holonomy, since it amounts to a simple variable switch in the
complete generating function.

A second useful operator is $\rev$, which switches the direction of a
step, $\rev(x,y)=(-x,-y)$. For example, $\rev(\NN)=\SS$, and
$\rev(\SW)=\NE$. The reverse of a set $\Y$, $\rev(\Y)$, is the set
resulting when $\rev$ is applied to each element of $\Y$. The reverse of
a walk $w=w_1w_2\dots$, is the application of $\rev$ to each step:
$\rev(w)=\rev(w_1)\rev(w_2)\dots$. This action can also be viewed as a
reflection in the line $x=-y$.

The reverse of a step set may result in a step set which has only the
trivial walk in the quarter plane. Aside from these, the effect of
this operation on the generating function can be determined, but as
it is much less direct than in the $\refl$ case, it is not entirely clear if it preserves
either holonomy or algebraicity. Evidence would seem to indicate
that in the case when the walk is non-trivial, that it does preserve
both properties. In Section~\ref{sec:future} we will give a better
intuition as to why this could be the case.

\subsection{Classifying step sets of cardinality three} In order to
classify all of the walks with step sets of cardinality three, it is
not necessary to consider all~$\binom83=56$ possibilities. Any step
set which is a subset of $\{\SE, \SS, \SW, \WW, \NW\}$ has no valid
walk in the quarter plane. These are 10 in total. Reflecting a step
set in the line $x=y$ yields a class of walks in bijection. There are
only four step sets which are invariant under this action, leaving 21
pairs of duplicates with respect to this action, and thus, there are
25 classes of walks, up to symmetry in the line $x=y$. We can reduce
this even further by determining other bijective classes. We shall
soon show that there are 11 non-bijective classes of walks over all,
and all of their generating functions can be classified. 

\subsection{Singular walks}\label{sec:singular}
\begin{figure} \center 
\includegraphics[width=1cm]{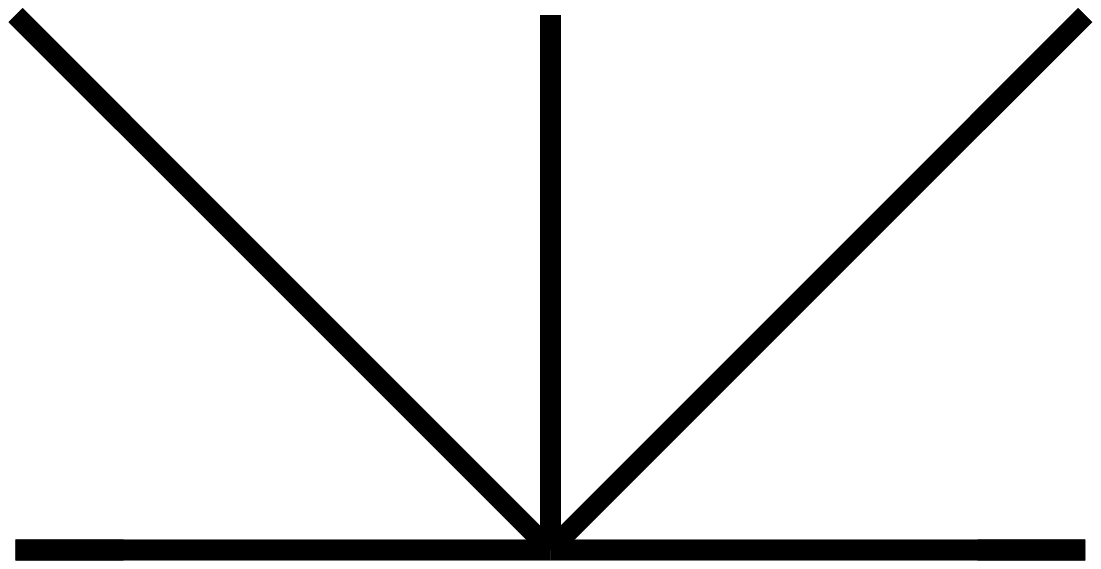}
\includegraphics[width=1cm]{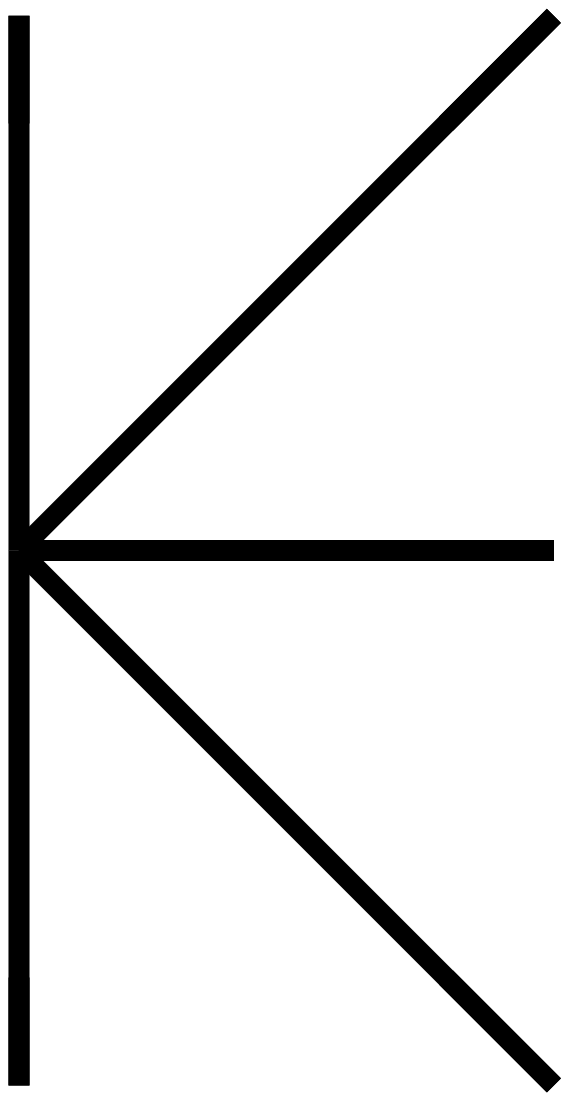} 
\includegraphics[width=1cm]{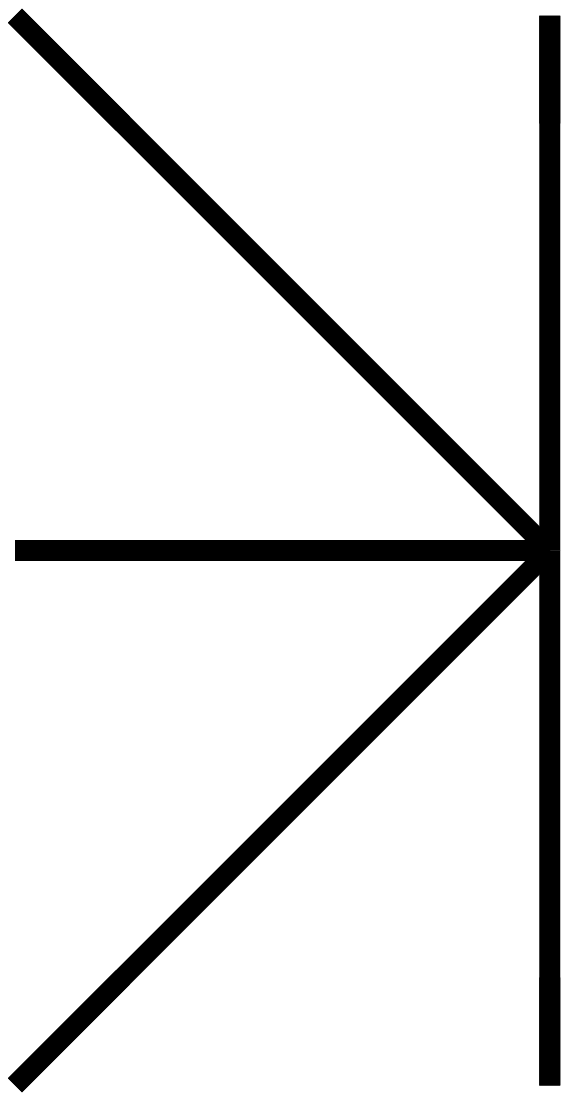}
\includegraphics[width=1cm]{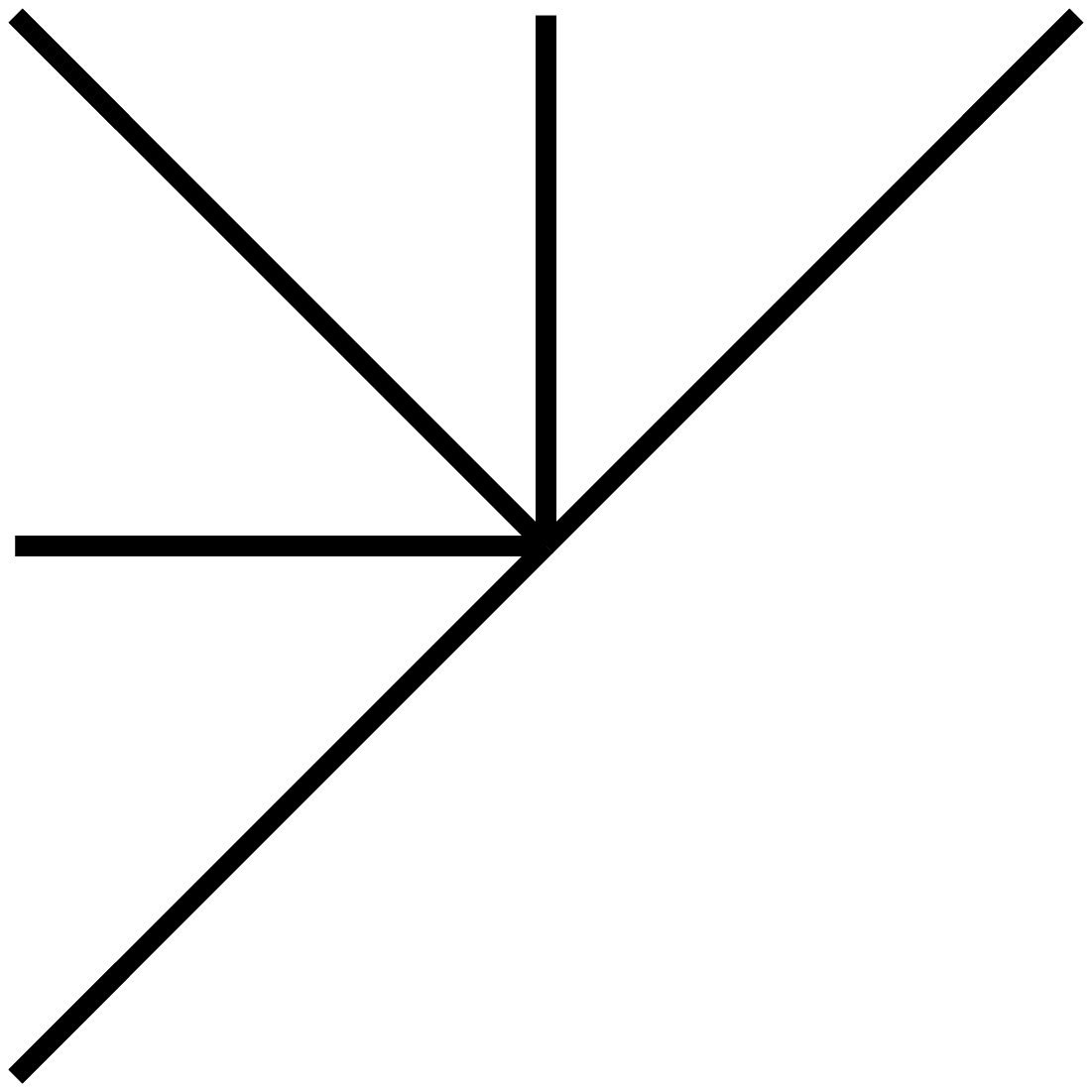}
\includegraphics[width=1cm]{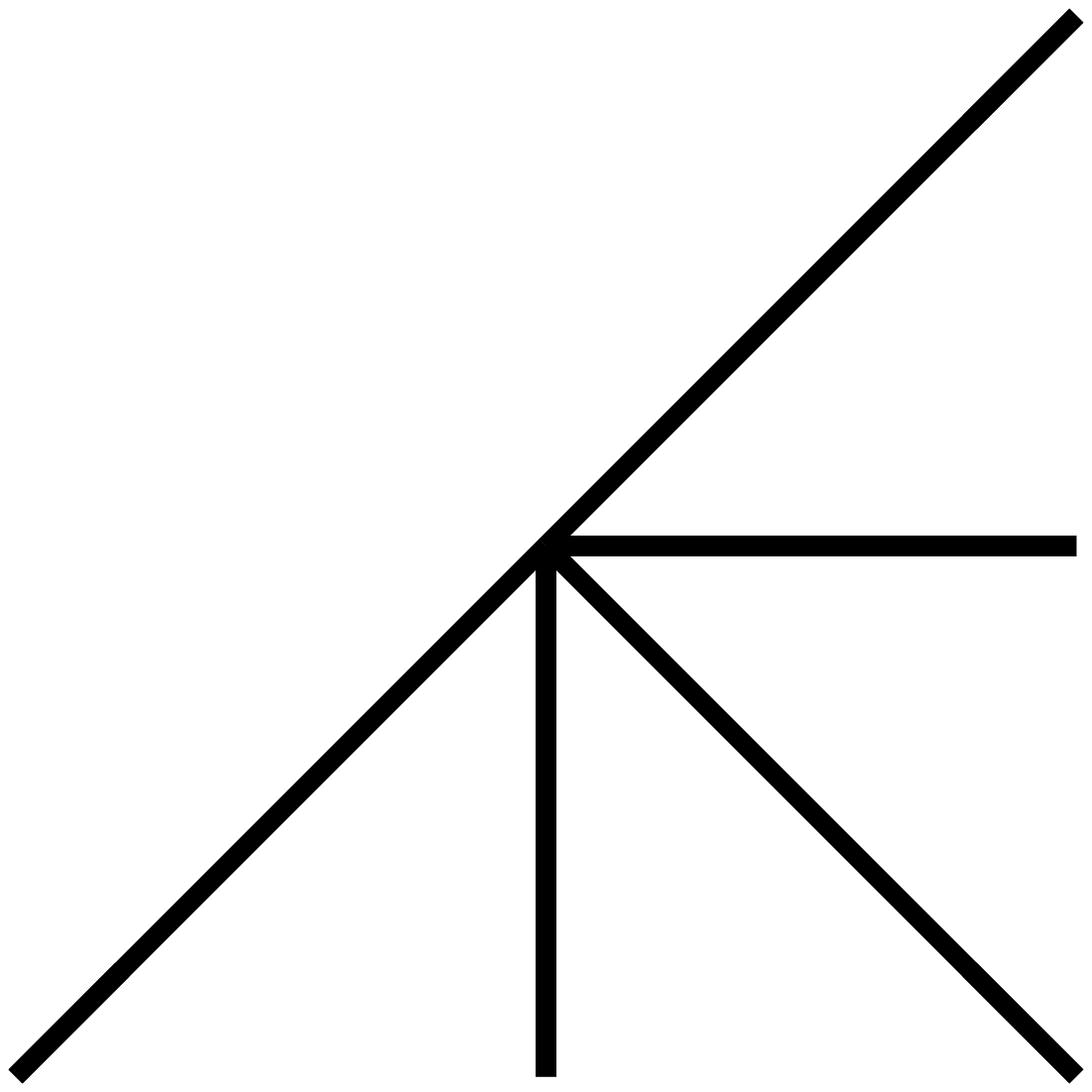}
\caption{Singular walks are subsets of these sets of directions} 
\label{fig:singular} \end{figure}

As we remarked in the introduction, each
step set is governed by at most two inequalities. In many cases one of
the inequalities implies the other, or one is trivial. For example,
the set $\{\NN, \NE, \EE\}$ only satisfies trivial inequalities, and
the vertical constraint on the walks~$\{\NN, \NE, \SS\}$ is implied by
the horizontal. These are both examples of {\em singular\/} step sets.
\begin{defn}[Singular step set] A step set is {\em singular} if it is
a subset of any of the following sets: \begin{enumerate} \item $\AA=\{\WW,
\NW, \NN, \NE, \EE\}$, $\refl(\AA)$, $\rev(\AA)$, $\refl(\rev(\AA))$;
\item $\BB=\{\NE, \NN, \NW, \WW, \SW \}$, $\refl(\BB)$.
\end{enumerate} These sets are pictured in vector format Figure~\ref{fig:singular}.
\end{defn} 

\begin{proposition}\label{thm:sing} If $\Y$ is a singular step set, then the complete
generating function $Q_\Y(x,y;t)$ is an algebraic function.
\end{proposition}

The algebraicity of the counting generating function for walks with a 
singular step set can be deduced from the observation that the
quarter plane condition of singular walks is equivalent to the half
plane condition, and these walks are algebraic by Theorem~\ref{thm:hpc}.
Furthermore, it is simple to construct a simple pushdown automaton
which recognizes the language of the walk. The language is thus unambiguously
context-free and it follows that the counting generating function is
algebraic~\cite{ChSc63}, and are easy to determine.

\begin{proof}
In order to prove the algebraicity of the complete generating function,
we give an explicit,
unambiguous grammar satisfied by the walks. A well-chosen weighting on
the certain variables gives a system of algebraic equations satisfied
by the complete generating function, and thus proves the algebraicity.

If the walk is singular, then there is one governing
inequality. Without loss of generality suppose it is the horizontal
condition. Divide the directions of the step set into three subsets:
Let $A=\{a_i\}_1^k=\{(j,1)\in\Y\}$, $B=\{b_i\}_1^l=\{(j,-1)\in\Y\}$,
$C=\{c_i\}_1^m=\{(j,0)\in\Y\}$. Thus, the vertical condition is given by
$\sum_{i=1}^k \# a_i \geq \sum_{j=1}^l\# b_j$.

If $\Y$ is singular, then~$\LY$ is generated by~$W$ in the
following grammar, which assures that the number of $A$s is always
greater than the number of $B$s for any prefix. The class $W$ is
decomposed by the last step up at a certain height. 
\[\begin{array}{rclrcl}
W &\rightarrow& (MA)^*M \qquad  & M &\rightarrow& \epsilon|CM|AMBM\\
A &\rightarrow& a_1|\dots|a_k   &  B &\rightarrow& b_1|\dots|b_l \\
C &\rightarrow& c_1|\dots|c_m\\
\end{array}\]
There is a direct correspondence between these grammars and the
functional equations satisfied by the generating function~\cite{ChSc63,
FlZiVa94}. Each step is weighted by a monomial corresponding to its
direction, namely $(i,j)$ is assigned to $x^iy^jt$. This gives a
solvable algebraic system for the complete generating system:
\begin{equation*}
\left\{
\begin{array}{rcl}
S(x,y,t)&=& \frac{M(x,y,t)^2A(x,y,t)}{1-M(x,y,t)A(x,y,t)};\\
M(x,y,t)&=& 1+ C(x,y,t)M(x,y,t)+A(x,y,t)B(x,y,t)M(x,y,t)^2;\\
A(x,y,t)&=& \sum_{i=1}^k a_i(x,y,t);
\quad B(x,y,t)= \sum_{i=1}^l b_i(x,y,t);\\
C(x,y,t)&=& \sum_{i=1}^m c_i(x,y,t).\\ 
\end{array}\right\}
\end{equation*}

We can then solve for~$Q_\Y(x,y;t)=S(x,y,t)$.
\end{proof}

To illustrate the process from the proof of
Proposition~\ref{thm:sing}, we determine algebraic equations satisfied
by the generating system for walks given by $\{\NE, \SW, \NN
\}$. Here, the horizontal constraint implies the vertical, and we set
$A=\{\NE\}$, $B=\{\SW\}$, and $C=\{\NN\}$. The walks are generated by
$S$ in the following grammar:
\[
S\rightarrow (MA)^*M \qquad M = \epsilon|CM|AMBM \]
In this case, the algebraic system is determined from the three
substitutions $A(x,y,t)= xyt$, $B(x,y,t)= \bar x \bar yt$, and
$C(x,y,t)= xt$ into the above system. The solution~$S(x,y;t)$ of this
system gives an expression for~$Q(x,y;t)$:
\[
Q(x,y;t)=S(x,y,t)=
-{\frac {-1+yt+\sqrt {1-2\,yt+{t}^{2}{y}^{2}-4\,{t}^{2}}}{t \left( 2\,
t-yx+{y}^{2}xt+yx\sqrt {1-2\,yt+{t}^{2}{y}^{2}-4\,{t}^{2}} \right) }}
.\]

In total there are 18 singular step sets. However, the walks which are
true two dimensional walks are governed by one of the following four
inequality types, and they are all isomorphic to one of four different
classes (Numbers refer to Table~\ref{tab:classes}):
\begin{equation*}
(1)\quad A, B, C\geq0, \qquad
(2)\quad A\geq B, C\geq 0, \qquad
(3)\quad  A+B\geq C, \qquad
(4)\quad A\geq B+C.
\end{equation*}

\begin{table}\small 
\begin{tabular}[t]{|c|c|c|c|c|}\hline
&$\Y$ & Counting GF& \begin{minipage}{4.5cm}\mbox{} \center Complete GF\\
$Q_k(x,y;t)=\sum a_{ij}(n) x^iy^jt^n$\\
\mbox{}\\
\end{minipage}
& cf.\\ \hline
1&\begin{minipage}{0.8cm}\includegraphics[height=.7cm]{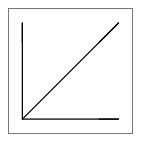}
\end{minipage}&
{$ \left( 1-3\,t \right)^{-1}$}
&
{$ \left( 1-t(x+y+xy) \right)^{-1}$}&
\S~\ref{sec:singular}
\\ \hline
%------------------------------------------------------------
2&\begin{minipage}{0.8cm}\includegraphics[height=.7cm]{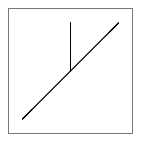}
\end{minipage}&
${\frac {-4\,t+1-\sqrt {-8\,{t}^{2}+1}}{4t(3t-1)}}$
& 
$-{\frac {-1+yt+\sqrt {1-2\,yt+{t}^{2}{y}^{2}-4\,{t}^{2}}}{t \left( 2\,
t-yx+{y}^{2}xt+yx\sqrt {1-2\,yt+{t}^{2}{y}^{2}-4\,{t}^{2}} \right) }}$
& \S~\ref{sec:singular}\\ \hline
%------------------------------------------------------------
3&\begin{minipage}{0.8cm}\includegraphics[height=.7cm]{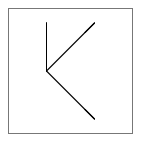}
\end{minipage}&
${\frac {-3\,t+1-\sqrt {-3\,{t}^{2}-2\,t+1}}{2t(3t-1)}}$
&
$-{\frac {-1+\sqrt {1-4\,x{t}^{2}-4\,{x}^{2}{t}^{2}}}{t \left( 2\,xt-y+
y\sqrt {1-4\,x{t}^{2}-4\,{x}^{2}{t}^{2}} \right)  \left( 1+x \right) }}$
&\S~\ref{sec:singular}\\ \hline
%------------------------------------------------------------
4&\begin{minipage}{0.8cm}\includegraphics[height=.7cm]{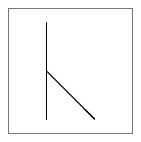}
\end{minipage}&
${\frac {-2\,t+1-\sqrt {-8\,{t}^{2}+1}}{2t(3t-1)}}$
&
$-{\frac {-1+\sqrt {1-4\,{t}^{2}-4\,x{t}^{2}}}{t \left( 2\,t+2\,xt-y+y
\sqrt {1-4\,{t}^{2}-4\,x{t}^{2}} \right) }}$
&\S~\ref{sec:singular}\\ \hline
%------------------------------------------------------------
5&\begin{minipage}{0.8cm}\includegraphics[height=.7cm]{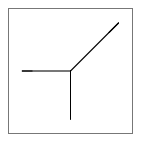}
\end{minipage}&
$\frac{T(1-t)+2t(T-1)\sqrt{1-T^2}}{tT(3t-1)}$
&
$\frac{xy-R(x,t)-R(y,t)}{xy-t(x+y+x^2y^2)}$
&\S~\ref{sec:krew}, ~\cite{MBM05}
\\ \hline
%------------------------------------------------------------
6&\begin{minipage}{0.8cm}\includegraphics[height=.7cm]{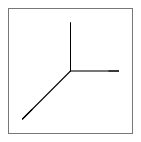}
\end{minipage}&
$\frac{\left(T^2t+T-2t \right)\sqrt{(1-T)(1+T^2/4+T^3/4)} +T +Tt}{tT(3t-1)}$

 & $\frac{xy-S(x,t)-S(y,t)}{xy-t(x^2y+xy^2+1)}$
&\S~\ref{sec:rkrew}
\\ \hline
%------------------------------------------------------------
7&\begin{minipage}{0.8cm}\includegraphics[height=.7cm]{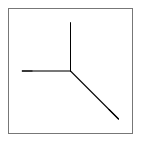}
\end{minipage}&
{$\frac{1-t-\sqrt{(1+t)(1-3t)}}{2t^2}$}& 
\begin{minipage}{5.2cm}\center
{\small $\frac{(i+1)(j+1)(i+j+2)n!}{\left(\frac{n-i-2j}{3}\right)!
\left(\frac{n-i+j}{3}+1\right)!\left(\frac{n+2i+j}{3}+2\right)!}
 $}
\end{minipage}
& \S~\ref{sec:florin}
\\ \hline
%------------------------------------------------------------
8&\begin{minipage}{0.8cm}\includegraphics[height=.7cm]{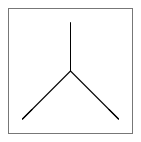}
\end{minipage}&
\multicolumn{2}{|c|}{Theorem~\ref{thm:DFIN1} (holonomic/transcendental)}
&\S~\ref{sec:nonsing} \\ \hline
%------------------------------------------------------------
9&\begin{minipage}{0.8cm}\includegraphics[height=.7cm]{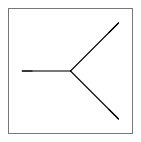}
\end{minipage}&
\multicolumn{2}{|c|}{Theorem~\ref{thm:DFIN2} (holonomic/transcendental)}
&\S~\ref{sec:nonsing}\\\hline
%-------------------------------------------------------------------
10&\begin{minipage}{0.8cm}\includegraphics[height=.7cm]{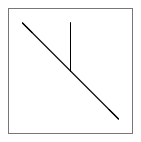}
\end{minipage}&
\multicolumn{2}{|c|}{Theorem~\ref{thm:notDfin} (Not holonomic)}
&\S~\ref{sec:notD}
\\ \hline
%------------------------------------------------------------
11&\begin{minipage}{0.8cm}\includegraphics[height=.7cm ]{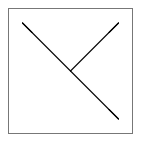}
\end{minipage}&
\multicolumn{2}{|c|}{Theorem~\ref{thm:Q11} (Not holonomic)}
&\S~\ref{sec:notD}
\\ \hline
%------------------------------------------------------------
\end{tabular} 

{\tiny
$T\equiv T(t)$ be the power series in $t$ defined by $T=t(2+T^3)$\\
$U=1-xT(1+T^3/4)+x^2T^2/4$\\
$R(x, t)=\frac1{2t}-\frac{1}{x}-\left(\frac1T-\frac1x\right)\sqrt{1-xT^2}$\\ 
$S(x,t)=\left( \frac{-2x}{T} \left( 1-{\frac {{T}^{2}}{2x}} \right)
+{\frac {1}{{tx}}} \right)\frac{\sqrt U}2 +\left( 1-tx-{\frac {t}{{x}^{2}}}
\right) \frac x{2t}$\\ }
\caption{Generating function data for walks whose step set is of cardinality three}
\label{tab:classes}
\end{table}

\subsection{Non-singular Walks}\label{sec:nonsing}
The remaining non-singular walks, (seven in total) break down as
follows. (Again, the numbers refer to those in Table 1.) 

Remark that step sets \#5 and \#6 are reverses of each other. Step set~
\#5 is known as Kreweras' walks after Kreweras' study~\cite{Kreweras65}. The
algebraicity of $W_5(t)$ is surprising and is well
studied. Step set~\#6 is examined in the next section, by applying the
same algebraic kernel method~\cite{MBM05} which gives the effective
algebraicity results for step set~\#5. The set~\#7 is examined in detail in
Section~\ref{sec:florin}. Step sets~\#8 and~\#9 satisfy the criteria
given in Theorem~\ref{thm:sym}, and are thus holonomic, and one can
exploit their symmetry in the kernel method to determine explicit
expressions for~$Q(x,y;t)$. As we shall see from calculations in
Section~\ref{sec:Df}, these walks are not algebraic.  The final two,
step sets~\#10 and~\#11, are not holonomic, and a proof is presented in Section~\ref{sec:notD}.

\section{The Algebraic Kernel Method} The technique we shall apply
here, in brief, uses the fundamental equation Eq.~\eqref{eqn:fund},
under different specializations of $x$ and $y$ which fix the coefficient
of $Q(x,y;t)$ in this equation, given by
$K_r:=1-t\sum_{(i,j)\in\Y}x^{i}y^{j}$.  The coefficient in this form
is called rational kernel, and is denoted $K_r(x,y)$. The
approach known as the {\em kernel method}, which has demonstrated a
growing popularity in combinatorics, see~\cite{MBM03, Prodinger03},
generally finds specializations of~$y$ as functions of~$x$ which
annihilate a polynomial form of the kernel (here, $xyK_r(x,y)$ would
be appropriate), and the problem is reduced to a simpler one. However,
Bousquet-M\'elou, in her study of Kreweras' walks~\cite{MBM05} takes a
slightly different approach, and instead finds a group of actions on
the pair $(x,y)$ which fixes the {\em rational kernel\/} $K_r(x,y)$.  This
group shall be known as the {\em group of the walk}.  These
specializations generate more functional equations from which we can
deduce information such as the algebraicity of $Q(x,y;t)$, or
potentially even explicit generating functions.

\subsection{The group of the walk}To each non-singular walk we associate a
group of actions which fix the rational kernel. 
For the walks we consider here, it will always be a
dihedral group, generated by two involutions. 2
We give an explicit definition for this particular case. 
\begin{defn}[The group of the walk ($G(\Y)$)]
Let $\Y$ be a fixed step set which is not singular. The {\em group of
  the walk \Y, denoted $G(\Y)$}, is the group of transformations which
map $R^2$ to itself generated by $\tau_x$ and $\tau_y$ defined as
follows.  Set
\[K(x,y):= xy K_r(x,y)=xy-t\sum_{(i,j)\in\Y}x^{i+1}y^{j+1}=a(x)ty^2+b(t,x)y+c(x)t,\]
and define $\tau_y$ as the transformation which maps a pair $(x,y)$ to
$(x, a(x)c(x)\bar y)$. Switch the roles of $x$ and $y$,
and likewise define $\tau_x$.
\end{defn}
\begin{table}\tiny
\center
\begin{tabular}[t]{|c|c|}\hline
%-------------------------------------------------
$\Y$ & Generators\\\hline
%-------------------------------------------------
\multicolumn{2}{|c|}{The Group is $D_2$}\\\hline
%-------------------------------------------------
% 1
\begin{minipage}{1.1cm}
\includegraphics[height=1.0cm]{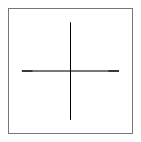}\end{minipage}
\begin{minipage}{1.1cm}
\includegraphics[height=1.0cm]{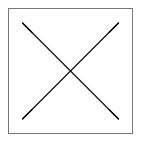}\end{minipage}
\begin{minipage}{1.1cm}
\includegraphics[height=1.0cm]{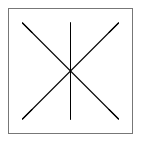}\end{minipage} 
\begin{minipage}{1.1cm}
\includegraphics[height=1.0cm]{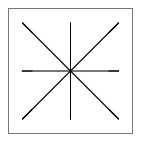}\end{minipage} &
 $\tau_x(x,y)=(\bar x,y),
\tau_y(x,y)=(x,\bar y)$\\ \hline

% 2
\begin{minipage}{1.1cm}
\includegraphics[height=1.0cm]{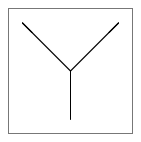}\end{minipage}
\begin{minipage}{1.1cm}
\includegraphics[height=1.0cm]{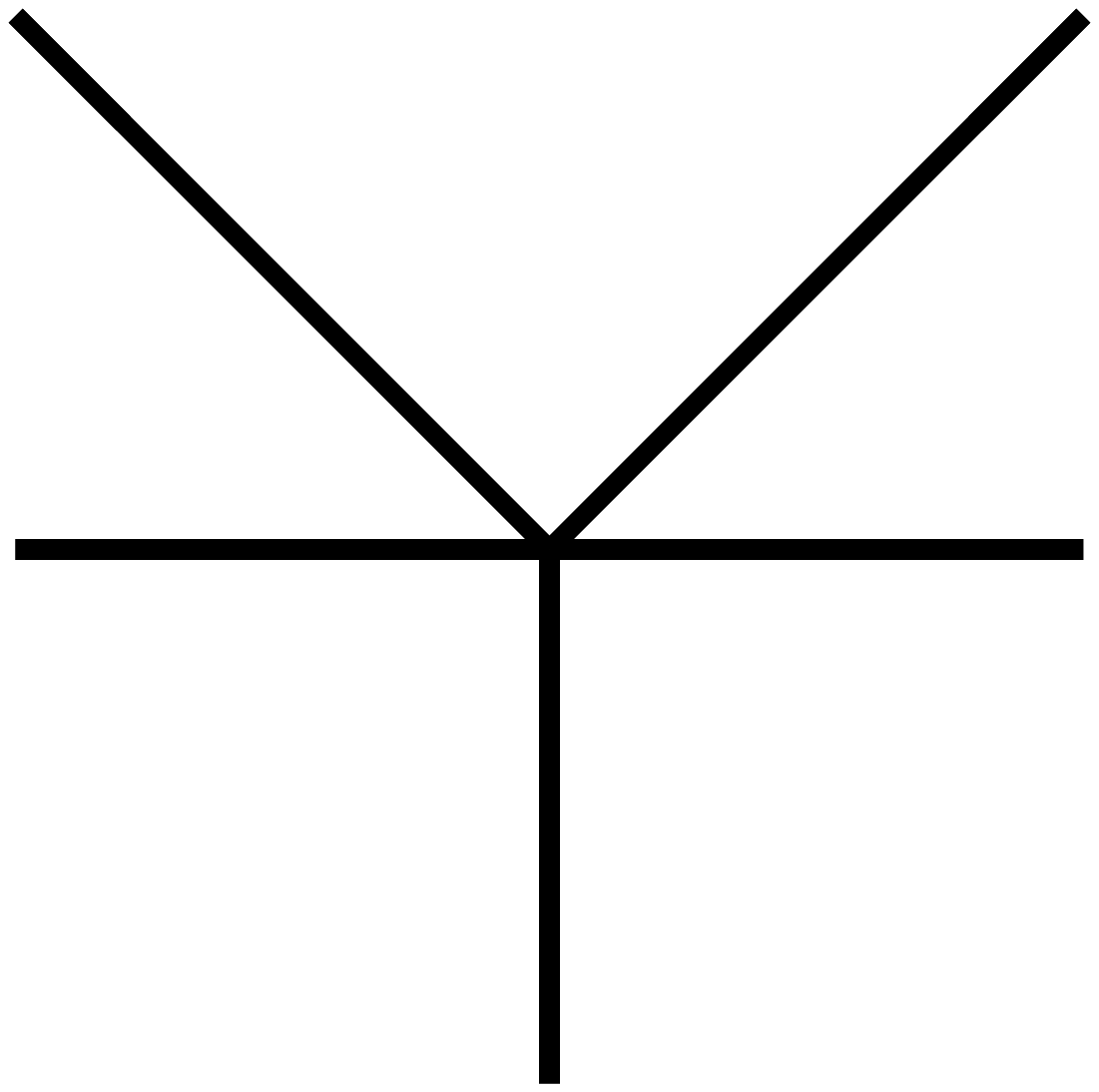}\end{minipage}&
 $\tau_x(x,y)=(\bar x,y),
\tau_y(x,y)=(x,{\frac {x}{y \left( {x}^{2}+1 \right) }})$ \\\hline

% 3
\begin{minipage}{1.1cm}
\includegraphics[height=1.0cm]{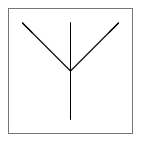}\end{minipage}
\begin{minipage}{1.1cm}
\includegraphics[height=1.0cm]{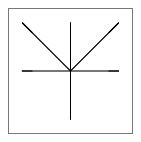}\end{minipage}&
 $\tau_x(x,y)=(\bar x,y),
\tau_y(x,y)=(x,{\frac {x}{y \left( {x}^{2}+1+x \right) }})$\\ \hline

% 4
\begin{minipage}{1.1cm}
\includegraphics[height=1.0cm]{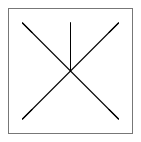}\end{minipage}
\begin{minipage}{1.1cm}
\includegraphics[height=1.0cm]{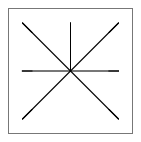}\end{minipage}&
 $\tau_x(x,y)=(\bar x,y),
\tau_y(x,y)=(x,{\frac {{x}^{2}+1}{y \left( {x}^{2}+1+x \right) }})$ \\ \hline

% 5
\begin{minipage}{1.1cm}
\includegraphics[height=1.0cm]{Plots/10010100}\end{minipage}
\begin{minipage}{1.1cm}
\includegraphics[height=1.0cm]{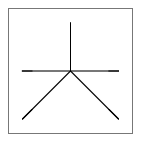}\end{minipage}&
 $\tau_x(x,y)=(\bar x,y),
\tau_y(x,y)=(x,{\frac {{x}^{2}+1}{xy}})$ \\ \hline

% 6
\begin{minipage}{1.1cm}
\includegraphics[height=1.0cm]{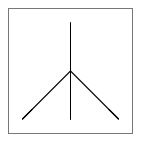}\end{minipage}
\begin{minipage}{1.1cm}
\includegraphics[height=1.0cm]{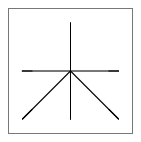}\end{minipage}&
 $\tau_x(x,y)=(\bar x,y),
\tau_y(x,y)=(x,{\frac {{x}^{2}+1+x}{xy}})$ \\ \hline

% 7
\begin{minipage}{1.1cm}
\includegraphics[height=1.0cm]{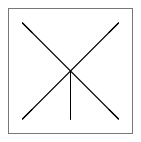}\end{minipage}
\begin{minipage}{1.1cm}
\includegraphics[height=1.0cm]{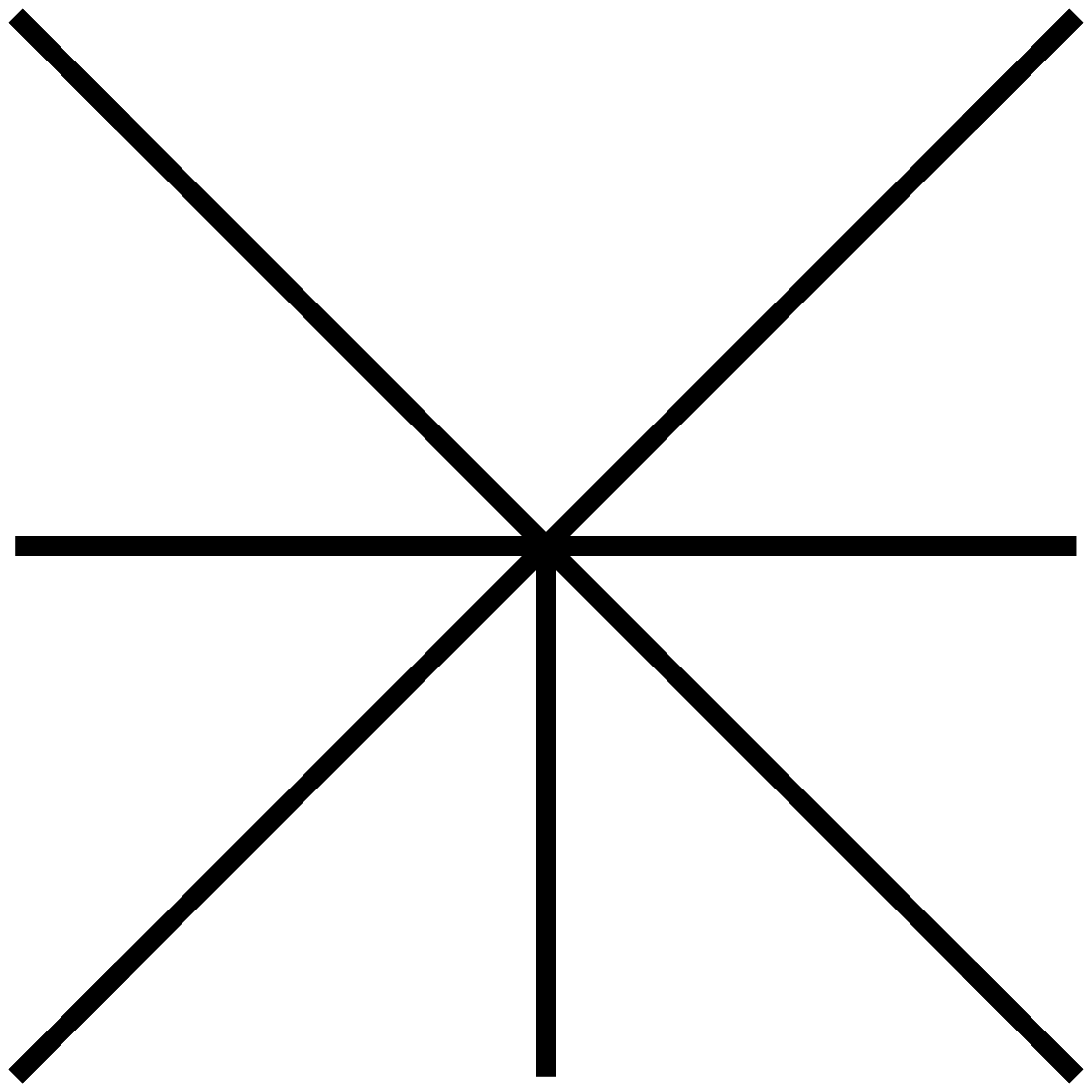}\end{minipage}&
$\tau_x(x,y)=(\bar x,y),
\tau_y(x,y)=(x,{\frac {{x}^{2}+1+x}{y \left( {x}^{2}+1 \right) }})$ \\ \hline

%-------------------------------------------------
\multicolumn{2}{|c|}{The Group is $D_3$}\\\hline
%-------------------------------------------------
% 8
\begin{minipage}{1.1cm}
\includegraphics[height=1.0cm]{Plots/10100100}\end{minipage}
\begin{minipage}{1.1cm}
\includegraphics[height=1.0cm]{Plots/01001010}\end{minipage}
\begin{minipage}{1.1cm}
\includegraphics[height=1.0cm]{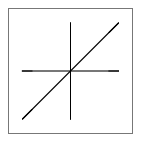}\end{minipage}
& 
$\tau_x(x,y)=({\frac {1}{xy}},y),
\tau_y(x,y)=(x,\frac {1}{xy})$\\ \hline

%9
\begin{minipage}{1.1cm}
\includegraphics[height=1.0cm]{Plots/10010010}\end{minipage}
\begin{minipage}{1.1cm}
\includegraphics[height=1.0cm]{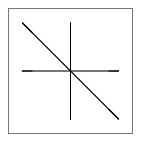}\end{minipage}&
 $\tau_x(x,y)=(y\bar x,y),\tau_y(x,y)=(x,x\bar y)$ \\ \hline
%-------------------------------------------------
\multicolumn{2}{|c|}{The Group is $D_4$}\\\hline
%-------------------------------------------------
% 10
\begin{minipage}{1.1cm}
\includegraphics[height=1.0cm]{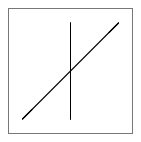}\end{minipage}&
 $\tau_x(x,y)=({\frac {1}{x{y}^{2}}},y),
\tau_y(x,y)=(x,{\frac {1}{xy}})$\\ \hline

% 11
\begin{minipage}{1.1cm}
\includegraphics[height=1.0cm]{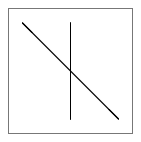}\end{minipage}&
 $\tau_x(x,y)=(\bar x y^2,y),
\tau_y(x,y)=(x, x\bar y)$ \\ \hline
\end{tabular}\\

\caption{The generators of the group of the step set in finite cases
  for all nearest neighbour walks.}
\label{tab:fingroup}
\end{table}

The generators for {\em all\/} step sets with finite groups (up to
symmetry in the line $x=y$) are given in Table~\ref{tab:fingroup}.
We can see that these generators fix the rational kernel:
 If $y=Y_0(x)$ and $y=Y_1(x)$ are the two
roots of the quadratic polynomial in $y$ given by $xyK_r(x,y)=0$, then $\tau_x(x,y)=(x, \frac{Y_1(x)Y_0(x)}y)$.
We can write $K_r=ta(x)\left(1-\frac{Y_0}{y}
\right)\biggl(y-Y_1\biggr)$, and thus
\begin{align*}
K_r(\tau_x(x,y))=K_r\left(x,\frac{Y_1Y_0}{y}\right) 
&=ta(x)\left(1-\frac{Y_0y}{Y_0Y_1}\right)\left(\frac{Y_0Y_1}{y}-Y_1\right)\\
&=ta(x)\left(1-\frac{y}{Y_1}\right)Y_1\left(\frac{Y_0}{y}-1\right)\\
&= K_r(x,y).
\end{align*}

In a later section we shall propose a conjecture for conditions on~$\Y$
which assure finiteness of the group generated by $\tau_x$ and $\tau_y$. 

This approach is inspired by the Galois automorphisms and the group of
the random walk describe by Fayolle {\em et al.\/} in~\cite{FaIaMa99},
and their treatment of the random walks. Here we shall apply this
method to step set~\#6, making use of some of the intermediary results
of Bousquet-M\'elou~\cite{MBM05}, from her study of step set \#5.

\subsection{Kreweras' Walks}\label{sec:krew}
The step set~$\Y=\{\NE, \SS, \WW\}$ is interesting since its
generating functions are algebraic, but~$\LY$ is {\em not\/}
context-free. (One can demonstrate this with the pumping lemma). As we
mentioned above, it has been studied in several different contexts,
and just recently a direct (i.e. combinatorial) explanation of its
algebraicity has been offered by Bernardi~\cite{Bernardi06}. He
demonstrates a bijection with a family of planar maps.
Here we shall stick to just
reporting enumerative data.  In the next section, however, we
follow the same method as Bousquet-M\'elou in~\cite{MBM05}, in order 
to determine a new family of walks with an algebraic generating function,
namely,~$\rev(\Y)$. It remains to be seen if the techniques of
Bernardi can be applied to $\rev(\Y)$ to explain algebraicity in that case. 

The fundamental equation expresses the complete generating function for
walks given by the set $\{\NE, \SS, \WW\}$ in terms of the walks which
returns to the $x$-axis. Theorem~1 of~\cite{MBM05} gives an explicit
formulation of this generating function. 
\begin{theorem}[Bousquet-M\'elou~\cite{MBM05}] 
\label{thm:kreweras}
Let $T\equiv T(t)$ be the power series
in $t$ defined by $T=t(2+T^3)$. The generating function for Kreweras'
walks ending on the $x$-axis is 
\[
Q_5(x,0,t)=\frac{1}{tx}\left(\frac1{2t}-\frac{1}{x}-\left(\frac1T-\frac1x\right)\sqrt{1-xT^2}\right).
\]
\end{theorem}
We have that $Q_5(y,0;t)= Q_5(0,y;t)$ and thus we have the following
expression for the complete generating function:
\[
Q_5(x,y;t)=\frac{xy-tQ_{5}(x,0;t)-tQ_{5}(0,y;t)}{xy-t(x+x+x^2y^2)}.
\]

\subsection{Reverse Kreweras walks}\label{sec:rkrew}
A second non-singular algebraic class is obtained from the reverse of
the Kreweras' walks. The language given by $\mathcal{L}(\Y)$ is
similarly not context-free. As Bousquet-M\'elou suggests~\cite{MBM-Labri}, the
algebraic kernel method can be applied in a straightforward manner to
this set. This results in explicit information about step set~\#6,
given by~$\{\NN, \EE, \SW \}$.

We have the following fundamental equation:
\[Q_6(x,y;t)=1+t\left(x+y+\bar x\bar y\right)Q_6(x,y;t)-\frac{t}{xy}\left(Q_6(0, y;t)+Q_6(x, 0; t)-Q_6(0,0; t)\right).  \]
We will now drop the index, rearrange the terms, and re-write this 
\begin{equation}\label{eq:fund}
xy K_rQ(x,y;t)=xy-t(R_0(x)+R_0(y)-R_{00}(t)),
\end{equation}
where the kernel of the fundamental equation is $K_r=(1-t(x+y+\bar
x\bar y))$, $R_0(x)=Q_6(x, 0; t)=Q_6(0, x; t)$ (by symmetry of the
step set), and $R_{00}(t)=Q_6(0,0;t)$.  First, we remark that since
the reverse Kreweras walks are reverses of Kreweras walks, the walks
of these two types that return to the origin are in bijection: simply
reverse the order of the path, and take the steps in the opposite
direction. Thus, $R_{00}(t)$ is a known algebraic power series which
we can deduce from Theorem~\ref{thm:kreweras}. The $x$ in the
denominator of the expression for $Q_5(x,0;t)$ is a removable
singularity, and hence we can determine $R_{00}=Q_5(0, 0;t)$ by taking
the limit as $x$ tends to 0. This gives the expression
\[ R_{00}(t)=\frac{4T-T^2}{8t}, \quad T=t(2+T^3).\]

Next, we determine the group of the walk. It is generated by
$\tau_y:(x,y)\mapsto(x, \frac1{xy})$ and $\tau_x:(x,y)\mapsto
(\frac1{xy}, y)$. We remark that this is the same as the group of the
walk for the Kreweras' walks.  Thus, we
have that $G(\Y)$ is $D_3$, the dihedral group on six elements.
We apply the invariance to obtain the following equalities
\begin{equation}\label{eq:ratk}
K_r(\bar x, \bar y)=K_r(\bar x, xy)=K_r(xy, \bar y),
\end{equation}
and note that~$K_r(\bar x, \bar y)$ is the rational kernel for
Kreweras' walks.  Denote this function by $\bar K_r\equiv K_r(\bar x,
\bar y)$.  Next, generate three equations substituting different
values for $x$ and $y$ (using Eq.~\eqref{eq:ratk}), into
Eq.~\eqref{eq:fund}:
\[
\begin{array}{ll} \bar x\bar y \bar K_r Q(\bar
x, \bar y; t)&=\bar x\bar y-t R_0(\bar x)-tR_0(\bar y)+ tR_{00}(t);\\ 
y\bar K_r Q(\bar x, xy; t)&=y-t R_0(\bar x)-tR_0(xy)+ tR_{00}(t);\\ 
x\bar K_r Q(xy,\bar y; t) &=x-t R_0(xy)-tR_0(\bar y)+ tR_{00}(t).
\end{array}
\] 
We form a composite equation taking the sum of the first two equations
and subtracting the third:
\begin{equation}\label{eq:comp} 
\bar x\bar y R(\bar x, \bar y)+yQ(\bar x, xy)-xQ(xy, \bar y)
=\frac{1}{\bar K_r} \left(\bar x\bar y+y -x -2tR_0(\bar x)+tR_{00}(t)\right).
\end{equation} 
The explicit expressions for $Y_0$ and $Y_1$ are 
\[ 
Y_0(x)=\frac{1-t\bar x-\sqrt{(1-t\bar x)^2-4t^2x}}{2tx}\qquad
\text{and}\qquad Y_1(x)
=\frac{1-t\bar x+\sqrt{(1-t\bar x)^2-4t^2x}}{2tx}.
\] 
Set~$\Delta$ to the common discriminant ($1-t\bar x)^2-4t^2x$ of $Y_0$
and $Y_1$.

The partial fraction expansion of $\bar K_r^{-1}$ is given by \[ \frac
1{\bar K_r}=\frac{1}{\left(1-Y_0\bar
y\right)\left(y-Y_1\right)}=\frac{1}{\sqrt{\Delta(x)}}
\left(\sum_{n\geq 0}\bar y^nY_0^n + \sum_{n\geq 1}y^nY_1^{-n} \right),
\] % 
The series $\Delta$ factors into three power series, respectively in
$C[x][\![t]\!], C[\![t]\!], C[\bar x][\![t]\!]$, which we shall call
$\Delta_+(x),\Delta(t)$, and $\Delta_-(\bar x)$. This is an instance
of a canonical factorization of a power series which is useful in many
enumeration problems. We now make direct use of some of
Bousquet-M\'elou's intermediary calculations. She determined that
$\Delta_{-}(\bar x)= 1-\bar x\left(T(1+T^3/4)+\bar xT^2/4 \right)$,
(which we later denote $1-\bar x \mathcal X$) with $T$ the unique
power series in $t$ defined by $T=t(2+T^3)$.

Next we extract the constant term with respect to $y$ from both sides
of Eq.~\eqref{eq:comp},
and express this using $Q_d(x)$, the generating series for walks that
end on the diagonal. This gives a new equation:
\begin{equation}
\label{eq:comp2}
-xQ_d(x)=\sqrt{\Delta(x)}^{-1}\biggl(2Y_0-x-2tR_0(\bar x)+tR_{00}(t)\biggr),
\end{equation}
Here, we have made use of the fact that $Y_1^{-1}\bar x=Y_0$. Now,
both the left and right hand sides are series in $\mathbb{C}[x, \bar
 x][[t]]$. For any such series $f(x, \bar x, t)$ denote the sub-series
whose coefficients of $t^n$ are polynomials in $\mathbb{C}[\bar x]$ by
$f^{\leq}$. We now isolate these sub-series from either side of Eq.~\eqref{eq:comp2}.
This gives
\begin{equation}\label{eqn:almost}
0=\frac{-2tR_0(\bar x)+tR_{00}(t)}{\sqrt{\Delta_0\Delta_{-}(\bar x)}}
- \left(\frac{x-2Y_0}{\sqrt{\Delta_0\Delta_{-}(\bar
x)}}\right)^{\leq}.
\end{equation}
We begin with the calculation 
\begin{eqnarray*}
\frac{x}{\sqrt{\Delta_0\Delta_{-}(\bar x)}}
&=&\frac{x}{\sqrt{\Delta_0}}\frac{1}{\sqrt{1-\bar
x\mathcal{X}}}\\
&=&\frac{x}{\sqrt{\Delta_0}}\left (1-\frac{\bar x\mathcal{X}}{2} +O(\bar x^2))\right)
\end{eqnarray*}
where $\mathcal{X}\in k[\![t]\!]$, leading to 
\begin{eqnarray*}
\left(\frac{x}{\sqrt{\Delta_0\Delta_{-}(\bar x)}}\right)^{\leq}
\biggl(\sqrt{\Delta_0\Delta_{-}(\bar x)}\biggr)
&=&\left(\frac{x}{\sqrt{\Delta_0\Delta_{-}(\bar x)}}-\frac{x}{\sqrt{\Delta_0}}\right)
\biggl(\sqrt{\Delta_0\Delta_{-}(\bar x)}\biggr)\\
&=&x\left(1-\sqrt{\Delta_{-}(\bar x)}\right).
\end{eqnarray*}
The remaining term is slightly more delicate to compute. We have
\begin{eqnarray*}
\left(\frac{2Y_0}{\sqrt{\Delta_0\Delta_{-}(\bar x)}}\right)^{\leq}
&=& \frac{\bar x(\frac{1-t^2}{t})}{\sqrt{\Delta_0\Delta_{-}(\bar x)}} -
    \left( \frac{\bar x}{t}\sqrt{\Delta_{+}}\right)^{\leq}.
\end{eqnarray*}
Next we develop $\frac{\bar x}{t}\sqrt{\Delta_{+}}$ to determine the
negative part
\begin{eqnarray*}
\left(\frac{\bar x}{t}\sqrt{\Delta_{+}(x)}\right)^{\leq}
&=&\left( \frac{\bar x}{t}\sqrt{1-xT^2} \right)^{\leq}\\
&=&\left(\frac{\bar x}t\left(1-xT^2/2-(xT^2)^2/4-O(x^6)\right)\right)^{\leq}\\
&=&\frac{\bar x}t-\frac{T^2}{2t}.
\end{eqnarray*}

Now we put it all together and we have 
\[\left(\frac{2Y_0}{\sqrt{\Delta_0\Delta_{-}(\bar x)}}\right)^{\leq}\sqrt{\Delta_0\Delta_{-}(\bar x)}=
\frac1t \left(\bar x(1-t^2)-(\bar x- \frac{T^2}{2})\sqrt{\Delta_0\Delta_{-}(\bar x)}\right)
\]

We now clear the denominator in Eq.~\eqref{eqn:almost} to get the
explicit expression given in the following theorem. We can indeed
conclude the algebraicity of $Q(x,0;t)$ from this expression. Using
this we reconstruct an expression for the complete
generating series demonstrating its algebraicity.
\begin{theorem} Let $\Y=\{\NN, \EE, \SW \}$. Then the generating series for walks that return to the x-axis is given by 
\[2Q_6(x,0;t)=Q_6(0,0;t)+
\left( \frac{-2x}{Tt} \left( 1-{\frac {{T}^{2}}{2x}} \right)
+{\frac {1}{{tx}}} \right)\sqrt {U} +\left( 1-tx-{\frac {t}{{x}^{2}}}
\right) x{t}^{-2} \] where~$T$ is the power series in~$t$ defined by
$T=t(2+T^3)$, and $U$ is the power series in $t$ defined by
$U(x)=1-xT(1+T^3/4)+x^2T^2/4$, and $Q_6(0,0;t)=(4T-T^2)/8t$.
\end{theorem}

\begin{theorem}
 Let $\Y=\{\NN, \EE, \SW \}$. Then the complete generating series for walks
with steps from $\Y$ is given by
\[Q_6(x,y;t)=\frac{xy-S(x,t)-S(y,t)}{xy-t(x^2y+xy^2+1)},\]
where $S(x,t)$ is given by $Q_6(x,0;t)$ from the previous theorem. 
\end{theorem}

\section{An interesting case}\label{sec:florin} 
\begin{figure}
\center \includegraphics[height=3cm]{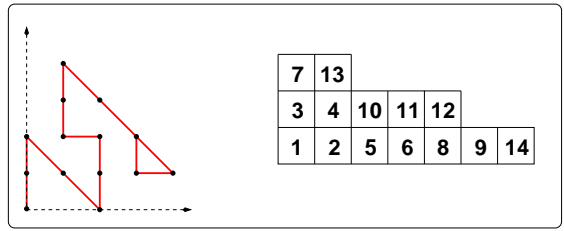}
\caption{A walk of step set~$\Y=\{\NN, \SE, \WW\}$ and its corresponding
Young tableau}
\label{fig:tandem}
\end{figure}

The walks given by
steps in $\Y=\{\NN, \SE, \WW\}$ are an interesting example because the
nature of the complete and the counting generating functions
differ. Regev~\cite{Regev81} proved that the number of such walks of
length~$n$ are equal to the $n^{th}$ Motzkin numbers, thus admitting
the counting generating function
\[W_7(t)=\frac{1-t-\sqrt{(1+t)(1-3t)}}{2t^2}.\] We can deduce this, and
an expression for the complete generating function, by exploiting a
bijective correspondence between these walks and standard Young
tableaux of height at most 3. Standard Young tableaux are labelled
Ferrer's diagrams of a partition such that the boxes are labelled in a
strictly increasing manner from left to right and from bottom to
top. We construct a standard Young tableau of height three of size $n$
from a walk $w=w_1, w_2, \ldots, w_n$ as follows. If $w_i=\NN$
(resp. \SE,\WW), then place label $i$ is in the next available space
to the right on the bottom row (resp. second, top). The final tableau
is increasing from left to right by construction, and the prefix
condition $\#\NN \geq \#\SE \geq \#\WW$ ensures that it is increasing
along the columns. Figure~\ref{fig:tandem} gives an example of such a
correspondence.

There is a well-known formula for counting standard Young tableaux,
known as the hook formula~\cite{Stanley99}. We apply this formula to
count the number of tableaux of
form $(n_1, n_2, n_3)$ from which we deduce the number of walks of
length $n$ that end at $(i,j)$.  The number of such tableaux are
\[
a(n_1, n_2, n_3)=(n_1-n_2+1)(n_2-n_3+1)(n_1-n_3+2)\frac{(n_1+n_2+n_3)!}{(n_1+2)!(n_2+1)!n_3!}.
\]
Now, the total length of the corresponding walk is the size of the
tableau, $n=n_1+n_2+n_3$ and to end at $(i,j)$, we have $i=n_2-n_3$
and $j=n_1-n_2$. Thus, when we make the substitution,
\[
a_{ij}(n)=\frac{(i+1)(j+1)(i+j+2)n!}{\left(\frac{n-i-2j}{3}\right)!\left(\frac{n-i+j+3}{3}\right)! \left(\frac{n+2i+j+6}{3}\right)!}.
\]

We can verify that this is {\em P-finite}~\cite{Lipshitz89}, and thus
the complete generating function is holonomic. 

\section{Explicit calculations for two holonomic classes}
\label{sec:Df}
In the next example we use the kernel method to determine enumerative
results for walks with steps from $\{\NN, \SE, \SW\}$ and $\{ \WW,
\NE, \SE \}$ . These results allow us to conclude that the series,
although holonomic, are transcendental. The property we exploit is
their axis symmetry.
Set $\Y=\{\NN,\SE,\SW\}$. The functional equation for $Q_\Y=Q$ is 
\begin{equation}\label{eq:FEQ}
K(x,y)Q(x,y;t)=xy-(x^2+1)tQ(x,0;t)-tQ(0,y;t)+tQ(0,0;t),
\end{equation}
with kernel $K(x,y)=xy(1-ty-t(x+\bar x)\bar y)$.  Let $Y_1(x)$ satisfy
$K(x,Y_1(x))=0$, and vanish at $t=0$. The quadratic formula gives
\begin{align*}
Y_1=Y_1(x)&=\frac{x-\sqrt{x^2-4x^3t^2-4xt^2}}{2xt}\\
&= \frac{1-\sqrt{1-4t^2(x+\bar x)}}{2t}=
(x+\bar x)t+(x+\bar x)^2t^3+O(t^5).
\end{align*} 
Next, we generate two equations based upon the observations that
$K(\bar x, y)=K(x,y)$ and $Y_1(\bar x)=Y_1(x)$:
\begin{align*}
0&=xY_1-(x^2+1)tQ(x,0;t)-tQ(0,Y_1;t)+tQ(0,0;t)\\
0&=\bar xY_1-(\bar x^2+1)tQ(\bar x,0;t)-tQ(0,Y_1;t)+tQ(0,0;t)
\end{align*}
The difference of these two equations yields a simpler third equation,
\begin{equation}\label{eq:FEQ2}
(x-\bar x)Y_1=tQ(x,0;t)(x^2+1)-tQ(\bar x,0;t)(\bar x^2+1).
\end{equation}

As we did in the case of the reverse Kreweras walks, we view both
sides of the equation as elements of $Q[x, \bar x][[t]]$, and isolate
the sub-series (the positive part, denoted by the $\geq$ in the
exponent) which is an element of $Q[x][[t]]$. We have that $Q(\bar x,
0;t)^\geq=Q(0,0;t)$, and thus the positive part of the right hand side
of Eq.~\eqref{eq:FEQ2} is $t(x^2+1)Q(x,0;t)-tQ(0,0;t)$. We shall
denote this series by~$H(x,t)$.  The positive
part of the left hand side requires a deeper series development. We
use the series expansion of $\sqrt{1-4t^2(x+\bar x)}$ to determine
the expression 
\[Y_1=\sum_{n>0}2(x+\bar x)^n\frac{t^{2n-1}}n\binom{2(n-1)}{(n-1)}.\]
The positive series $\left((x+\bar
x)^n\right)^\geq=\sum_{k=0}^{\lfloor\frac n2\rfloor}\binom nk
x^{n-2k}$. 

Thus, when we combine these two, we can describe $a(n,k)$, the coefficient of
$t^{2n-1}x^k$ in $\left((x-\bar x)Y_1\right)^\geq$:
\[
a_n=0, k\equiv n \mod 2; \binom{n+1}{(n+k+1)/2}\frac{k}{n+1},
\]
where $C_n$ is the $n^\text{th}$ Catalan number $C_n=\binom {2n}n\frac{1}{n+1}$.
We substitute this result into the orginal equation. 
\begin{equation}\label{eq:FEQ3}
K(x,y)Q(x,y;t)=xy-H(x,t)-tQ(0,y;t).
\end{equation}
Next, we define $X_1$ as the root $K(X_1, y)$ which  vanishes at
$t=0$. Again, this is easily determined using the quadratic
formula. Substituting this into Eq.~\eqref{eq:FEQ3}, we have
that $tQ(0,y)=X_1y-H(X_1,t)$.  This gives the following theorem.

\begin{theorem}\label{thm:DFIN1}
Let $\Y=\{\NN, \SE, \SW\}$. Then the complete generating series for walks with steps
from $\Y$ is given by
\[
Q_8(x,y;t)=\frac{xy-H(x,t)-M(y,t)+H(M(y,t),t)}{xy-txy^2-t(x^2+1)}
\]
where $H(x,t)$ is the power series given by
$[x^kt^{2n-1}]H(x,t)=\binom{n+1}{(n+k+1)/2}\frac{k}{n+1}$ when $k\equiv n \pmod 2$ and 0 otherwise; and $M(y,t)=\frac1{2t}{y-t{y}^{2}-\sqrt {{y}^{2}-2\,{y}^{3}t+{t}^{2}{y}^{4}-4\,{
t}^{2}}}
=\bar yt+t^2+(y+\bar y^3)t^3+(y^2+3\bar y^2)t^4+O \left( {t}^{5} \right) )
$. 
\end{theorem}

A similar calculation gives the following theorem. 
\begin{theorem}\label{thm:DFIN2}
Let $\Y=\{\NE, \SE, \WW\}$. Then the generating series for walks with steps
from $\Y$ returning to the $y$-axis is given by
\[
R(y,t)=Q_9(0,y;t)= \frac{(y^2-1)S(y,t)}{y^2+1}
\]
where the non-zero coefficients of $S(y,t)$ are given by
$[y^{n-2k}t^{2n-1}]S(y,t)=2\binom nk \binom {2n-2}{n-1}\frac 1n$, for
integer $0\leq k\leq n/2$.

The complete generating series is given by 
\[
Q_9(x,y;t)=\frac{xy-T(x-R(T,t))-tyR(y,t)}{xy-tx^2(y^2+1)-ty}
\]
where $T=\frac{(x-t)-\sqrt{(x-t)^2-4t^2x^4}}{tx^2}$.  
\end{theorem}

\subsection{A word about transcendence}
Suppose, for the purpose of illustrating a contradiction, that
$Q_8(x,y;t)$ were algebraic. Were it so, this implies the algebraicity
of $Q_8(x,0;t)$ and $Q_8(0,0;t)$, which then in turn implies the
algebraicity of $H(1,t)$ from the statement of the theorem.  However, if $H(1,t)=\sum
a(n) t^{2n-1}$, we can show that asymptotically, as $n$ tends to
infinity $a(n)$ tends to $\frac{8^n\sqrt{2}}{4\pi n^2}$. Thus by
Criterion~D of \cite{Flajolet87}, the series $H(1,t)$ is
transcendental. Thus, we have established a contradiction, and
$Q_8(x,y;t)$ is not algebraic. Recall, however, that it is holonomic,
by Theorem~\ref{thm:sym}.

\section{Walks which are not holonomic}\label{sec:notD}
Knowing that walks in the half plane and slit planes are algebraic,
and considering that walks in the quarter plane are ultimately
relatively simple objects-- they are governed by only two
inequalities, it may be surprising that there are classes with
generating functions which are not holonomic. However, it is already
known that Bousquet-M\'elou and Petkov\v sec's knight's
walks~\cite{BoPe03} ($\Y=\{(2,-1), (-1, 2)\}$) are not holonomic, and,
in fact, evidence would seem to indicate that many walks in the
quarter plane are not holonomic.

The final two walks we consider have non-holonomic complete generating
functions. This non-holonomy is established in~\cite{MiReXX}, and is
based on a similar problem of self-avoiding walks in wedges~\cite{JaPrRe06}. The
general argument is to show that there are an infinite number of
singularities in the counting generating function. This is similar to
the argument of Bousquet-M\'elou and Petkov\v sec, whose proof also
use an adaptation of a kernel method argument. However, for
these walks, it seems that there is a more direct application of the
group of the walk. Specifically, we use the fact that the group of the
walk is infinite to demonstrate the infinite set of singularities.

From~\cite{MiReXX} we pull the following two results which completes
our classification. In the next section we present a small summary of
the argument. 

\begin{theorem}\label{thm:notDfin}
Neither the complete generating function~$Q_{10}(x,y;t)$, nor the counting
generating function $W_{10}(t)=Q_{10}(1,1;t)$ of nearest-neighbour
walks in the first quadrant with steps from $\{\NE,\SE,\NW\}$ 
are holonomic functions with respect to their variable sets.
\end{theorem}

\begin{theorem}\label{thm:Q11}
Neither the complete generating function~$Q_{11}(x,y;t)$, nor the counting
generating function $W_{11}(t)=Q_{11}(1,1;t)$ of nearest-neighbour
walks in the first quadrant with steps from $\{\NE,\SE,\NN\}$ 
are holonomic functions with respect to their variable sets.
\end{theorem}

\subsection{The iterated kernel method and the step set
  $\Y=\{\NE,\SE,\NW\}$} To begin, we recall the fundamental equation of these walks, 
\begin{align*}
  Q_{10}(x,y;t) &= 1 + t xy Q_{10}(x,y;t) 
  + t {x}{\bar y} \left(Q_{10}(x,y;t) - Q_{10}(x,0;t) \right)\\
  &+ t {y}{\bar x} \left(Q_{10}(x,y;t) - Q_{10}(0,y;t) \right),
\end{align*}
and its kernel form,
\begin{equation}\label{eqn:kform}
 \left(xy - tx^2y^2 - tx^2 - ty^2\right) Q(x,y) = xy - tx^2 Q(x,0) - ty^2Q(y,0).
\end{equation}
Here, for brevity we did write $Q(x,y;t)$ as $Q(x,y)$, and have
used the $x \leftrightarrow y$ symmetry to rewrite $Q(0,y)$ as
$Q(y,0)$. There are two solutions for the
kernel~$K(x,y)=xy - tx^2y^2 - tx^2 - ty^2$ as a function of $y$
\begin{equation}
  \y_{\pm 1}(x) = \frac{x}{2t(1+x^2)} \left(1 \mp \sqrt{1-4t^2(1+x^2)} \right).
\end{equation}
Since we can show that
\begin{equation*}
  \y_{+1} ( \y_{-1} (x) ) = x\quad \text{ and }\quad \y_{-1} ( \y_{+1} (x) ) = x,
\end{equation*}
if we  write $\y_n(x) =  \left( \y_1 \circ \right)^n(x)$, and likewise
for $\y_{-n}$, we have that the set $\{\y_n|n\in\Z\}$ forms an infinite group,
under the operation $\y_n\circ \y_m=\y_{n+m}$, with identity~$\y_0=x$. 
We can also show the useful relation,
 \begin{equation*}
   \frac{1}{\y_1(x)} + \frac{1}{\y_{-1}(x)} = \frac{1}{tx},
 \end{equation*}
and its iterated equivalent, 
 \begin{equation}\label{eqn:recip}
   \frac{1}{\y_n(x)} = \frac{1}{t \y_{n-1} } - \frac{1}{\y_{n-2}(x)}.
 \end{equation}

By construction we have $K(x,\y_1(x)) = 0$, thus substituting $y =
\y_1(x)$ into Eq.~\eqref{eqn:kform} gives (after a little tidying)
\begin{equation*}
  Q(x,0) = \frac{\y_1(x)}{x} \frac{1}{t} - \frac{\y_1(x)^2}{x^2} Q(\y_1(x), 0).
\end{equation*}
Now we substitute~$x = \y_{n}(x)$ into this equation to obtain
\begin{equation*}
  Q(\y_n(x),0) = \left(\frac{\y_{n+1}(x)}{\y_n(x)}\right) \frac{1}{t} 
  - \left(\frac{\y_{n+1}(x)}{\y_{n}(x)}\right)^2 Q(\y_{n+1}(x), 0). 
\end{equation*}
Using this expression for $Q(\y_n(x),0)$ for various $n$, we can
iteratively generate a new expression for $Q(x,0)$:
\begin{equation*}
  Q(x,0) = \frac{1}{x^2 t} \sum_{n=0}^{N-1} (-1)^n\y_n(x) \y_{n+1}(x)
  + (-1)^N \left(\frac{\y_N(x)}{x} \right)^2 Q(\y_N(x),0).
\end{equation*}
Since $\y_n(x) = xt^n + o(xt^n)$ we have that $\ds \lim_{N \to
  \infty} \y_N(x) = 0$ as a formal power series in $t$ and
consequently
\begin{equation}\label{eqn:Qxo}
  Q(x,0) = \frac{1}{x^2 t} \sum_{n \geq 0} (-1)^n \y_n(x) \y_{n+1}(x).
\end{equation}
We now specialize the fundamental equation:
\begin{equation}\label{eq:infsing}
\left(1-3t\right)  Q(1,1;t) = 1-2tQ(1,0;t)
=1 - 2\sum_{n \geq 0} (-1)^n \y_n(1) \y_{n+1}(1).
\end{equation}

The result follows from the following lemma, which is proved in~\cite{MiReXX}
\begin{lemma}\label{thm:work}
Suppose $q_c$ is a zero of
$\bar\y_N(q):=\y_N(1;\frac{q}{1+q^2})^{-1}$, and that $q_c\neq 0$. Then 
\begin{enumerate}
\item For all $k \neq N$,  $\bar\y_k(q_c) \neq 0$; and 
\item The function $Q(1,0;\frac{q}{1+q^2})$ has a pole at the $q = q_c$.
\end{enumerate}
\end{lemma}

We now have all the components in place to prove the main result.

\begin{proof}[Proof of Theorem~\ref{thm:notDfin}]
The function $Q(1,1; \frac{q}{1+q^2})$
has a set of poles given by the zeroes of the
$\bar\y_n(q)$, by the preceding lemma. The set of such poles form an
infinite set. Thus, $Q(1,1, \frac{q}{1+q^2})$ is not holonomic. For a multivariate series to be holonomic, any of its algebraic
specializations must be holonomic, and as  $Q(1,1, \frac{q}{1+q^2})$ is
an algebraic specialization of both $Q(x,y;t)$ and $Q(1,1;t)$, neither
of these two functions are holonomic either. 
\end{proof}

A natural extension of this work would be to complete the
story, and use the functional equations to find enumerative data.

\subsubsection{What happens with the method when the generating function
  IS holonomic?}
This is a natural question, and it speaks to robustness of the
method. We can in general follow the steps to express the counting
generating function as a sum of products of $\y_n$. However, when the
group of the walk is finite, the $\y_n$ form a finite group, and thus
this sum does not converge as power series. 

The difficulty in using this method arises in showing that the
singularities do not cancel.

\section{Prospects for general criteria}\label{sec:future}
The group of the walk is useful for determining additional equations
in the the kernel method. However, an examination of all nearest
neighbour walks and their groups reveals a surprising fact about when
this group is finite. Table~\ref{tab:fingroup} lists all nearest
neighbour walks with finite groups (up to symmetry in the line $x=y$) 
\begin{proposition}
Consider {nearest neighbour} walks in the quarter plane,  and further
suppose \Y~is not singular. Then $G(\Y)$ is finite if and only if one
of the following is true:
\begin{enumerate}
\item $\Y$ is $x$-- or $y$-- axis symmetric;
\item $\Y= \rev(\Y)$;
\item $\Y = \refl(\rev(\Y))$
\item $\Y$ is either Kreweras or $\rev$(Kreweras)
\end{enumerate}
\end{proposition}
The following conjecture is true in the case of $|\Y|=3$, and all
other known cases. 
\begin{conjecture} Suppose the step set group \Y~is not singular.
Then  $Q_\Y$ is holonomic if and only if the group of its step set is finite.
\end{conjecture}
In fact, in view of the details behind the iterated kernel method which we used to prove the
non-holonomy of the two final step sets, this seems to be a
reasonable conjecture indeed. 

Combining the above proposition and conjecture, we have potential
conditions on $\Y$ for it to have a holonomic complete generating function.
\begin{conjecture}
%If $Q_\Y$ is transcendental, then
The generating function $Q_\Y$ is holonomic if and only if at
least one of the following holds:
\begin{itemize}
\item $\Y$ is singular;
\item $\Y$ is $x$-- or $y$-- axis symmetric;
\item $\Y=\rev(\Y)$ (path reversibility);
\item $\Y = \refl(\rev(\Y))$.
\end{itemize}
\end{conjecture}

This conjecture also implies that $\rev$ preserves holonomy. 
If $K_r(x,y)$ is the kernel related to $\Y$, then $K_r(\bar x,
\bar y)$ is the kernel related to $\rev(\Y)$. By taking the viewpoint
that $\bar x$ and $\bar y$ are mere variables, we see that if the
\Y-step set group is finite, then so is the $\rev(\Y)$-step set as the
groups will be the same.

We could equally conjecture similar results for $W_\Y(t)$.

We note that step set \#10, whose generating function is not
holonomic, does have some symmetry-- $\Y=\refl(\Y)$, but this is not
sufficient to offer holonomy. It is worth noting that the other
two classes with this same symmetry are Kreweras and reverse
Kreweras.

If true, this conjecture would respond positively to a conjecture of
Gessel on the nature of the walks given by the step set~$\Y=\{\NN,
\SE, \SS, \NW\}$.

\subsection{Future Work}
We would like to characterize the classification
combinatorially. Several natural generalizations are possible, such as what is a meaningful
definition for group in 3 dimensions? Does an analogous definition
still yield dihedral groups? What about step sets which are not of the
nearest neighbour type; what kinds of groups come out of that? We are
also currently investigating a similar classification for different
wedge shapes, such as 1/8-plane and 3/4-plane. Fayolle et al. give a
characterization for when the group is finite. Can this be translated
into straightforward combinatorial terms?

\subsection{Acknowledgments} This work was completed, in part, during an NSERC
(Natural Sciences and Engineering Research Council of Canada) Postdoc
at LABRI, Universit\'e Bordeaux I, under the guidance of Mireille
Bousquet-M\'elou.

% \bibliographystyle{plain}
% \bibliography{Main}

\end{document}